# Parameter estimation of ODE's via nonparametric estimators

**Nicolas J-B. Brunel**[*]

*Laboratoire IBISC, Université d'Evry*
*523, place des Terrasses, 91025 EVRY cedex, France.*
*Ecole Nationale Supérieure d'Informatique pour l'Industrie et l'Entreprise*
*1, square de la Résistance, 91025 EVRY cedex, France.*
*e-mail:* `nicolas.brunel@ibisc.univ-evry.fr`

**Abstract:** Ordinary differential equations (ODE's) are widespread models in physics, chemistry and biology. In particular, this mathematical formalism is used for describing the evolution of complex systems and it might consist of high-dimensional sets of coupled nonlinear differential equations. In this setting, we propose a general method for estimating the parameters indexing ODE's from times series. Our method is able to alleviate the computational difficulties encountered by the classical parametric methods. These difficulties are due to the implicit definition of the model. We propose the use of a nonparametric estimator of regression functions as a first-step in the construction of an M-estimator, and we show the consistency of the derived estimator under general conditions. In the case of spline estimators, we prove asymptotic normality, and that the rate of convergence is the usual $\sqrt{n}$-rate for parametric estimators. Some perspectives of refinements of this new family of parametric estimators are given.

## 1. Introduction

Ordinary differential equations are used for the modelling of dynamic processes in physics, engineering, chemistry, biology, and so forth. In particular, such a formalism is used for the description of ecosystems (for example competing species in biology), or of cell regulatory systems such as signaling pathways and gene regulatory networks [12]. Usually, the model for the state variables $x = (x_1, \ldots, x_d)^\top$ consists of an initial value problem

$$\begin{cases} \dot{x}(t) & = F(t, x(t), \theta), \ \forall t \in [0,1], \\ x(0) & = x_0, \end{cases} \quad (1.1)$$

[*]This paper was written during a visit at EURANDOM, Technical University of Eindhoven and a postdoctoral fellowship at IBISC, Université d'Evry. It was partly supported by Genopole®, the European PASCAL network, Université d'Evry and EURANDOM. Part of this work was presented at the workshop "Parameter Estimation in Systems Biology", Glasgow, March 2007.

where $F$ is a time-dependent vector field from $\mathbb{R}^d$ to $\mathbb{R}^d$, $d \in \mathbb{N}$, and $\theta \in \Theta$, $\Theta$ being a subset of a Euclidean space. When data are available such as a time series, we are interested in the problem of estimation of the coefficients parametrizing the ODE. In principle, this may be done by some classical parametric estimators, usually the least squares estimator [28] or the Maximum Likelihood estimator (MLE). Different estimators have been derived in order to take into account some particular features of the differential equation such as special boundary values (there exists a function $g$ linking the values at the boundary i.e. $g(x(0), x(1)) = 0$ instead of the simple initial value problem), or random initial values or random parameters [9]. Otherwise, there may be some variations on the observational process such as noisy observation times that necessitate the introduction of appropriate minimization criteria [26].

Despite their satisfactory theoretical properties, the efficiency of these estimators may be dramatically degraded in practice by computational problems that arise from the implicit and nonlinear definition of the model. Indeed, these estimators give rise to nonlinear optimization problems that necessitate the approximation of the solution of the ODE and the exploration of the (usually high-dimensional) parameter space. Hence, we have to face possibly numerous local optima and a huge computation time. Instead of considering the estimation of $\theta$ straightforwardly as a parametric problem, it may be useful to look at it as the estimation of a univariate regression function $t \mapsto x(t)$ that belongs to the (finite dimensional) family of functions satisfying (1.1). Alternative approaches to MLE has been used in applications, such as two-step methods: in a first step, a proxy for the solution of the ODE is obtained by nonparametric methods, and in a second step, estimates of the ODE parameter are derived by minimizing a given distance. Varah [40] initiated this approach by using and differentiating a smooth approximation of the solution based on least-squares splines as [38], or Madar *et al.* [29] with cubic splines (and a well-chosen sequence of knots). In the same spirit, approximation of the solutions of the ODE are provided by smoothing methods such as local polynomial regression [22, 30, 11], or neural networks [42], within the same two-step estimation scheme. Slight modifications of these approaches are the adaptation of collocation methods to statistical estimation where the solution is approximated by Lagrange polynomials [20]. This two-step approach has also been considered in the functional data analysis (FDA) framework proposed Ramsay and Silverman [34], which is based on the transformation of data into functions with smoothing cubic splines. Ramsay proposed Principal differential Analysis (PDA) [32] for the estimation of linear ODE's. PDA has been recently extended to nonlinear ODE's and ameliorated by repeating the two steps, introducing iterated PDA (iPDA) [31, 41]. Recently, the iPDA approach has then been extended to a so-called generalized smoothing approach by Ramsay *et al.*, [33]. In this algorithm, the smoothing and the estimation of the parameter ODE are considered jointly, and Ramsay *et al.* proposed an estimation algorithm inspired from profile likelihood methods. The paper provides a versatile estimator, as a lucid account of the current status and open questions concerning the statistical estimation of differential equations. We refer to [33] and the subsequent discussions for a broad view of the difficulties for parameter estimation of ODE's. The use of nonparametric estimators is motivated by the computational simplicity of the estimation method, but an additional motivation is that the functional point of view enables one to use prior knowledge on the solutions of the ODE such as positivity or boundedness whereas it is difficult to exploit the strictly parametric form. Indeed, it implies that we have a thorough knowledge of the influence of the parameters on the qualitative behavior of the solutions of (1.1), which is rarely the case. In this paper, we study the consistency and asymptotic properties of this general two-step estimator based on nonparametric estimators of the trajectories of the ODE. We focus in particular on the use of spline regression in the first step, as it is one of the estimator used previously in the construction of the proxy, but we discuss also the generality and potential ameliorations of this approach.

In the next section, we introduce the statistical model and we define the so-called two-step estimator of $\theta$. We show that under broad conditions this estimator is consistent, and we discuss some straightforward extensions of this estimator to different cases. In section 3, we derive an asymptotic representation of the two-step estimator, and we focus on the case of the least squares splines. We discuss also the generality of this result with respect to other types of nonparametric estimators. In the last section, we give some simulation results obtained with the classical Lotka-Volterra's population model coming from biology.

## 2. Two-step estimator

### 2.1. Statistical model

We want to estimate the parameter $\theta$ of the ordinary differential equation (1.1) from noisy observations at $n$ points in $[0,1]$, $0 \leq t_1 < \cdots < t_n \leq 1$,

$$y_i = x(t_i) + \epsilon_i, \; i = 1, \ldots, n, \tag{2.1}$$

where the $\epsilon_i$'s are i.i.d centered random variables. The ODE is indexed by a parameter $\theta \in \Theta \subset \mathbb{R}^p$ and initial value $x_0$; the true parameter value is $\theta^*$ and the corresponding solution of (1.1) is $x^*$.

The vector field defining the ODE is a function $F : [0,1] \times \mathcal{X} \times \Theta \to \mathbb{R}^d$ ($\mathcal{X} \subset \mathbb{R}^d$) of class $C^m$ w.r.t $x$ for every $\theta$ and with $m \geq 1$. It is a Lipschitz function so that we have existence and uniqueness of a solution $x_{\theta,x_0}$ to (1.1) on a neighborhood of 0 for each $\theta$ and $x_0$; and we assume that every corresponding solution can be defined on $[0,1]$. Hence, the solutions $x_{\theta,x_0}$ belong to $\mathcal{C}^{m+1}([0,1], \mathbb{R}^d)$. Moreover, we suppose also that $F$ is a smooth function in $\theta$ so that each solution $x_{\theta,x_0}$ depends smoothly[1] on the parameters $\theta$ and $x_0$. Then, we suppose that $F$ is of class $C^1$ in $\theta$ for every $x$. Let $f_\Sigma$ be the density of the

---
[1] if $F$ depends smoothly on $x$ and $\theta$ then the solution depends on the parameter by the same order of smoothness, see Anosov & Arnold, Dynamical systems, p.17.

noise $\epsilon$, then the log-likelihood in the i.i.d case is

$$l(\theta, x_0, \Sigma) = \sum_{i=1}^{n} \log f_\Sigma(y_i - x_{\theta,x_0}(t_i)) \qquad (2.2)$$

and the model that we want to identify is parametrized by $(\theta, x_0, \Sigma) \in \Theta \times \mathcal{X} \times \mathcal{S}^+$ for instance when the noise is centered Gaussian with covariance matrix $\Sigma$ ($\mathcal{S}^+$ is the set of symmetric positive matrices). An alternative parametrization is $(\theta, x_{\theta,x_0}, \Sigma) \in \Theta \times \mathcal{F} \times \mathcal{S}^+$, with $\mathcal{F}$ the set of functions that solve (1.1) for some $\theta$ and $x_0$, thanks to the injective mapping between initial conditions and a solution.

In most applications, we are not really interested in the initial conditions but rather in the parameter $\theta$, so that $x_0$ or $x_{\theta,x_0}$ can be viewed as a nuisance parameter like the covariance matrix $\Sigma$ of the noise. We want to define estimators of the "true" parameters $(x^*, \theta^*)$ ($x^* = x_{\theta^*, x_0^*}$) that will be denoted by $(\hat{x}_n, \hat{\theta}_n)$. The estimation problem appears as a standard parametric problem that can be dealt with by the classical theory in order to provide good estimators (with good properties, e.g. $\sqrt{n}$-consistency) such as the Maximum Likelihood Estimator (MLE). Indeed, from the smoothness properties of $F$, the log-likelihood $l(\theta, x_0)$ is at least $C^1$ w.r.t $(\theta, x_0)$ so that we can define the score $s(\theta, x_0) = (\frac{\partial l}{\partial \theta}^\top \frac{\partial l}{\partial x_0}^\top)^\top$. If $s(\theta, x_0)$ is square integrable under the true probability $P_{(x^*, \theta^*)}$, we can claim under weak conditions (e.g. theorem 5.39 [39]) that the MLE is an asymptotically efficient estimator. The difficulty of this approach is then essentially practical because of the implicit dependence of $x$ on the parameter $(\theta, x_0)$, which prohibits proper maximization of $l(\theta, x_0)$. Indeed, derivative-based methods like Newton-Raphson are not easy to handle then and evaluation of the likelihood necessitates the integration of the ODE, which becomes a burden when we have to explore a huge parameter space. Moreover, the ODE's proposed for modelling may be expected to give a particular qualitative behavior which can be easily interpreted in terms of systems theory, e.g. convergence to an equilibrium state or oscillations. Typically, these qualitative properties of ODE are hard to control and involve bifurcation analysis [25] and may necessitate a mathematical knowledge which is not always accessible for huge systems. Moreover, boundedness of the solution $x^*$ ($a \leq x^*(t) \leq b$, with $a, b \in \mathbb{R}^d$) may be difficult to use during the estimation via the classical device of a constraint optimization. Hence, these remarks motivate us to consider the estimation of an ODE as a functional estimation and use flexible methods coming from nonparametric regression from which we could derive a likely parameter for the ODE.

## 2.2. Principle

We use consistent nonparametric estimators of the solution $x^*$ and its derivative $\dot{x}^*$ in order to derive a fitting criterion for the ODE and subsequently the M-estimator of $\theta^*$ corresponding to the criterion. We denote by $\|f\|_{q,w} =$ $\left(\int_0^1 |f(t)|^q w(t)dt\right)^{1/q}, 0 < q \leq \infty$, the $L^q(w)$ norm on the space of integrable functions on $[0,1]$ w.r.t. the measure $w$ (and $L^q$ is the classical norm with respect to Lebesgue measure). We suppose that $w$ is a continuous and positive function on $[0,1]$. The principle of two-step estimators is motivated by the fact that it is rather easy to construct consistent estimators of $x^*$ and of its derivative $\dot{x}^*$. The estimation of the regression function and its derivatives is a largely addressed subject and several types of estimators can be used such as smoothing splines (more generally smoothing in Reproducing Kernel Hilbert Spaces [43, 3]), kernel and local polynomial regression [13], series estimators [10]. So, one can construct a consistent estimator $\hat{x}_n$ of $x^*$ that can achieve the optimal rate of convergence with an appropriate choice of regularization parameters. One can derive also from $\hat{x}_n$ an estimator of the derivative $\dot{x}^*$ by differentiating directly the smoothed curved $\hat{x}_n$ so that we have $\hat{\dot{x}}_n = \dot{\hat{x}}_n$. This simple device provides consistent estimator of the derivative for Nadaraya-Watson estimators [17], spline regression [45], smoothing splines [43], but other consistent estimators of the derivative can be obtained from local polynomial regression or wavelet estimator (e.g [6]). In this section, we consider then simply that we have two general nonparametric estimators $\hat{x}_n$ and $\hat{\dot{x}}_n$, such that $\|\hat{x}_n - x^*\|_q = o_P(1)$ and $\|\hat{\dot{x}}_n - \dot{x}^*\|_q = o_P(1)$. Hence, for the estimation of the parameter $\theta$, we may choose as criterion to minimize the function $R_{n,w}^q(\theta) = \|\hat{\dot{x}}_n - F(t, \hat{x}_n, \theta)\|_{q,w}$ from which we derive the two-step estimator

$$\hat{\theta}_n = \arg\min_\theta R_{n,w}^q(\theta). \tag{2.3}$$

Thanks to the previous convergence results and under additional suitable conditions to be specified below, we can show that

$$R_{n,w}^q(\theta) \to R_w^q(\theta) = \|\dot{x}^* - F(\cdot, x^*, \theta)\|_{q,w}$$

in probability, and that this discrepancy measure enables us to construct a consistent estimator $\hat{\theta}_n$.

We are left with three choices of practical and theoretical importance: the choice of $q$, the choice of $w$ and the choice of the nonparametric estimators. In this paper, we show the consistency in a general situation, but we provide a detailed analysis of the asymptotics of $\hat{\theta}_n$ when $q = 2$, $\hat{x}_n$ is a regression spline (least-squares splines, with a number of knots depending on the number of observations $n$), and when $\hat{\dot{x}}_n = \dot{\hat{x}}_n$. In that case, the two-step estimation procedure is then computationally attractive as we have to deal with two (relatively simple) least-squares optimization problems (linear regression in the first step and nonlinear regression in the second step). As a consequence, this is a common approach in applications and we will put emphasis on the importance of the weight function $w$ in asymptotics. Despite the limitation to regression splines, it is likely that our result can be extended to classical nonparametric estimators, as the ones listed above.

### 2.3. Consistency

We show that the minimization of $R_{n,w}^q(\theta)$ gives a consistent estimator for $\theta$. We introduce the asymptotic criterion

$$R_w^q(\theta) = \left( \int_0^1 (F(\cdot, x^*, \theta^*) - F(\cdot, x^*, \theta))^q \, w \, dt \right)^{1/q}$$

derived from $R_{n,w}^q$ and we make the additional assumption:

$$\forall \epsilon > 0, \inf_{\|\theta - \theta^*\| \geq \epsilon} R_w^q(\theta) > R_w^q(\theta^*), \tag{2.4}$$

which may be viewed as an identifiability criterion for the model.

**Proposition 2.1.** *We suppose there exists a compact set $\mathcal{K} \subset \mathcal{X}$ such that $\forall \theta \in \Theta, \forall x_0 \in \mathcal{X}, \forall t \in [0,1]$, $x_{\theta,x_0}(t)$ is in $\mathcal{K}$, and $w$ is a positive continuous function on $[0,1]$. Moreover we suppose that uniformly in $(t,\theta) \in [0,1] \times \Theta$, $F(t,\cdot,\theta)$ is $K-$ Lipschitz on $\mathcal{K}$. If $\hat{x}_n$ and $\hat{\dot{x}}_n$ are consistent, and $\hat{x}_n(t) \in \mathcal{K}$ almost surely, then we have*

$$\sup_{\theta \in \Theta} |R_{n,w}^q(\theta) - R_w^q(\theta)| = o_P(1).$$

*Moreover, if the identifiability condition (2.4) is fulfilled the two-step estimator is consistent, i.e.*

$$\hat{\theta}_n - \theta^* = o_P(1).$$

*Proof.* In order to show the convergence of $|R_{n,w}^q(\theta) - R_w^q(\theta)| = |\|\hat{\dot{x}}_n - F(\cdot, \hat{x}_n, \theta)\|_{w,q} - \|F(\cdot, x^*, \theta) - F(\cdot, x^*, \theta^*)\|_q|$, we make the following decomposition

$$\begin{aligned} |R_{n,w}^q(\theta) - R_w^q(\theta)| &\leq \| \left( \hat{\dot{x}}_n - F(\cdot, \hat{x}_n, \theta) \right) + (F(\cdot, x^*, \theta) - F(\cdot, x^*, \theta^*)) \|_{q,w} \\ &\leq \|\hat{\dot{x}}_n - F(\cdot, x^*, \theta^*)\|_{q,w} + \|F(\cdot, \hat{x}_n, \theta) - F(\cdot, x^*, \theta)\|_{q,w}. \end{aligned} \tag{2.5}$$

Since all the solutions $x_{\theta,x_0}(t)$ and $\hat{x}_n(t)$ stay in $\mathcal{K} \subset \mathcal{X}$, and $x \mapsto F(t,x,\theta)$ are $K-$ Lipschitz uniformly in $(t,\theta)$, we obtain for all $(t,\theta) \in [0,1] \times \Theta$

$$\|F(\cdot, \hat{x}_n, \theta) - F(\cdot, x^*, \theta)\|_{q,w} \leq KM \left( \int_0^1 |\hat{x}_n(t) - x^*(t)|^q dt \right)^{1/q} = KM \|\hat{x}_n - x^*\|_q. \tag{2.6}$$

where $M$ is a upper bound for $w$. Together, (2.5) and (2.6) imply

$$\begin{aligned} \sup_{\theta \in \Theta} |R_{n,w}^q(\theta) - R_w^q(\theta)| &\leq \|\hat{\dot{x}}_n - F(\cdot, x^*, \theta^*)\|_{w,q} \\ &\quad + \sup_{\theta \in \Theta} \|F(\cdot, \hat{x}_n, \theta) - F(\cdot, x^*, \theta)\|_{w,q} \\ &\leq M \times \|\hat{\dot{x}}_n - F(x^*, \theta^*)\|_q + KM \times \|\hat{x}_n - x^*\|_q. \end{aligned}$$

and consequently, by the consistency of $\hat{x}_n$ and $\dot{\hat{x}}_n$,

$$\sup_{\theta \in \Theta} |R^q_{n,w}(\theta) - R^q_w(\theta)| = o_P(1).$$

With the additional identifiability condition (2.4) for the vector field $F$, Theorem 5.7 of [39] implies that the estimator $\hat{\theta}_n$ converges in probability to $\theta^*$. □

The minimization of a distance between a nonparametric estimator and a parametric family is a classical device for parametric estimation (Minimum Distance Estimators, MDE). For instance, it has been used for density estimation with Hellinger distance [2], or for regression with $L^2$ distance in the framework of regression model checking [24]. The difference between two-estimators and MDE's is that two-step estimators do not minimize directly the distance between the regression function and a parametric model, but between the derivative of the regression function and the parametric model. This slight difference causes some differences concerning the asymptotic behavior of the parametric estimators. Nevertheless, two-step estimators are closer to MDE's than the generalized smoothing approach of Ramsay *et al.* [33] that uses also parametric and nonparametric estimators. These ingredients are used in a different manner, as the criterion maximized by the generalized smoothing estimator can be seen as the Lagrangian relaxation of a constrained optimization problem in a space of cubic splines (with given knots): the parametric and nonparametric estimators solve

$$\min_{g,\theta} \sum_{i=1}^n |y_i - g(t_i)|_2^2 \text{ subject to } \|\dot{g} - F(\cdot, g, \theta)\|_2^2 < \epsilon.$$

The smoothing approach finds a curve that solves approximately a differential equation, which is then close to the spline approximation of the solution of the ODE (with parameter $\hat{\theta}$) computed by collocation [7]. Moreover, the optimization problem is solved iteratively which implies that the nonparametric estimator is computed adaptively with respect to the parametric model and the data are repeatedly during estimation, whereas two-step estimators used the data once without exploiting the parametric model.

**Remark 2.1** (**Local identifiability**). The global identifiability condition (2.4) is difficult to check in practice and in all generality. In the case $q = 2$, we can focus on a simpler criterion by considering the Hessian $J^*$ of $R^2_w(\theta)$ at $\theta = \theta^*$:

$$J^* = \int_0^1 (D_2 F(t, x^*(t), \theta^*))^\top D_2 F(t, x^*(t), \theta^*) w(t) dt.$$

If $J^*$ is then nonsingular, the asymptotic criterion behaves like a positive definite quadratic form on a neighborhood $\mathcal{V}(\theta^*)$ of $\theta^*$, so the criterion $R^2$ separate the parameters. It is then critical to use a positive weight function $w$, or else some part of the path $t \mapsto x^*(t)$ might be useless for the identification of the identification and discrimination of the parameter. At the opposite, $w$ might be used to emphasize some differences between parameters. As it is shown in the next section, the value of $w$ at the boundaries has a dramatic influence on the asymptotics of the two-step estimator.

**Remark 2.2** (**Random times and errors in times**). If the observation times $t_1, \ldots, t_n$ are realizations of i.i.d. random variables $(T_1, \ldots, T_n)$ with common c.d.f $Q$, the nonparametric estimators $\hat{x}_n$, as the one used before, are relevant candidates for the definition of the two-step estimator since they are still consistent under some additional assumptions on $Q$.

As in the setting considered by Lalam and Klaassen [26], the observation times may be observed with some random errors $\tau_i = t_i + \eta_i, i = 1, \ldots, n$, (the $\eta_i$'s being some white noise) so we have to estimate $x$ from the noisy bivariate measurements $(\tau_i, y_i)$. Consistent nonparametric estimators have been proposed for the so-called "errors-in-variables" regression and some examples are kernel estimators [14] and splines estimators [23] (in the $L^2$ sense). Hence, we can define exactly the same criterion function $R_n^2$ and derive a consistent parametric estimator.

### 2.4. Partially observed ODE

The estimator proposed can be easily extended to cases where several variables are not observed. Indeed, if the differential system (1.1) is partially linear in the following sense

$$\begin{cases} \dot{u} &= G(u, v; \eta) \\ \dot{v} &= H(u; \eta) + Av \end{cases} \quad (2.7)$$

with $x = (u^\top v^\top)^\top$, $u \in \mathbb{R}^{d_1}$, being observed, $v \in \mathbb{R}^{d_2}$ being unobserved, and $d_1 + d_2 = d$ (the initial conditions are $x_0 = (u_0^\top v_0^\top)^\top$), i.e. $x(t_i)$ is replaced by $u(t_i)$ in (2.1) (the noise $\epsilon_i$ being then $d_1$-dimensional). We want to estimate the parameter $\theta = (\eta, A)$ when $H$ is a nonlinear function and $A$ is a matrix, so we can take advantage of the linearity in $v$ in order to derive an estimator for $v$. We can derive a nonparametric estimator for $v$ by using $\hat{u}_n$ and the fact that $t \mapsto v(t)$ must be the solution of the non-homogeneous linear ODE $\dot{v} = Av + H(\hat{u}_n; \eta)$, $v(0) = v_0$. The solution of this ODE is given by Duhamel's formula [19]

$$\forall t \in [0, 1], \ \hat{v}_n(t) = \exp(tA) v_0 + \int_0^t \exp((t-s)A) H(\hat{u}_n(s); \theta) ds, \quad (2.8)$$

which then can be plugged into the criterion $R_n^q(\theta)$. This estimator depends explicitly on the initial condition $v_0$ which must be estimated at the same time

$$(\hat{\theta}, \hat{v}_0) = \arg \min_{(\eta, A, v_0)} R_n^q(\eta, A, v_0) = \left\| \dot{\hat{u}}_n - F(\hat{u}_n, \hat{v}_n; \eta) \right\|_{w,q}.$$

As previously, if $H$ is uniformly Lipschitz the integral $\int_0^t \exp((t-s)A) H(\hat{u}_n; \theta) ds$ converges (uniformly) in probability in the $L^q$ sense to $\int_0^t \exp((t-s)A) H(u^*; \theta) ds$ as soon as $\hat{u}_n$ does, hence $R_n^q(\theta, A, v_0)$ converges also uniformly to the asymptotic criterion

$$R_w^q(\theta, v_0) = \| \dot{u}^* - F(u^*, v^*; \theta) \|_{w,q}.$$

The estimator $(\hat{\theta}, \hat{v}_0)$ is consistent as soon as $R_w^q(\theta, v_0)$ satisfies the identifiability criterion (2.4).

## 3. Asymptotics of two-step estimators

We give in this section an analysis of the asymptotics of two-step estimators. We focus on the least squares criterion $R_{n,w}^2$ when the estimator of the derivative is $\hat{\dot{x}}_n$. In that case, we show that the two-step estimator behaves as the sum of two linear functionals of $\hat{x}_n$ of different nature: a smooth and a non-smooth one. Under broad conditions, we show that the two-step estimator is asymptotically normal and that we can derive the rate of convergence. In particular, we show that it is the case for the regression splines. A key result is that the use of a weight function $w$ with boundary conditions $w(0) = w(1) = 0$ makes the non-smooth part vanish, implying that the two-step estimator can have a parametric rate of convergence. In the contrary, the two-step estimator has a nonparametric rate.

### 3.1. Asymptotic representation of two-step estimators

We derive an asymptotic representation of $\hat{\theta}_n$ by linearization techniques based on Taylor's expansion of $R_{n,w}^2$ around $\theta^*$. This shows that the two-step estimator behaves as a simple linear functional of $\hat{x}_n$. We introduce the differentials of $F$ at $(x, \theta)$ w.r.t. $x$ and $\theta$ and we denote them $D_1 F(x, \theta)$ and $D_2 F(x, \theta)$ respectively. For short, we adopt the notation $D_i F^*$, $i = 1, 2$ for the functions $t \mapsto D_i F(x^*(t), \theta^*)$, $i = 1, 2$, and $D_{12} F^*$ for $t \mapsto D_1 D_2 F(x^*(t), \theta^*)$.

**Proposition 3.1.** *We suppose that $D_{12}$ exists and $w$ is differentiable. We introduce the two linear operators $\Gamma_{s,w}$ and $\Gamma_{b,w}$ defined by*

$$\Gamma_{s,w}(x) = \int_0^1 \left( D_2 F^{*\top}(t) D_1 F^*(t) w(t) - \frac{d}{dt}(D_2 F^*(t) w(t)) \right) x(t) dt$$

*and*

$$\Gamma_{b,w}(x) = w(1) D_2 F^{*\top}(1) x(1) - w(0) D_2 F^{*\top}(0) x(0)$$

*If $D_1 F$, $D_2 F$ are Lipschitz in $(x, \theta)$, $J^*$ is invertible, and $\hat{x}_n, \hat{\dot{x}}_n$ are (resp.) consistent estimators of $x^*$ and $\dot{x}^*$, then*

$$\hat{\theta}_n - \theta^* = J^{*-1} \left( \Gamma_{s,w}(\hat{x}_n) - \Gamma_{s,w}(x^*) + \Gamma_{b,w}(\hat{x}_n) - \Gamma_{b,w}(x^*) \right) + o_P(1).$$

*Proof.* For sake of simplicity, we suppose that the vector field $F$ does not depend on $t$, but the proof remains unchanged in the non-autonomous case. We remove also the dependence in $t$ and $n$ for notational convenience and introduce $F^*$ and $F(\hat{x}, \theta^*)$. The first order condition $\nabla_\theta R_n^2(\hat{\theta}) = 0$ implies that $\int_0^1 (D_2 F(\hat{x}, \hat{\theta}))^\top (\hat{\dot{x}} - F(\hat{x}(t), \hat{\theta})) w \, dt = 0$. This gives

$$\int_0^1 \left( D_2 F(\hat{x}, \hat{\theta}) \right)^\top \left( (\hat{\dot{x}} - \dot{x}^*) + F^* - F(\hat{x}, \theta^*) + F(\hat{x}, \theta^*) - F(\hat{x}, \hat{\theta}) \right) w \, dt = 0$$

and

$$\int_0^1 \left(D_2F(\hat{x},\hat{\theta})\right)^\top \left((\dot{\hat{x}}-\dot{x}^*) + D_1F(\tilde{x}^*,\theta^*)(x^*-\hat{x}) + D_2F(\hat{x},\tilde{\theta}^*)(\theta^*-\hat{\theta})\right) w\, dt = 0$$

with $\tilde{x}^*$ and $\tilde{\theta}^*$ being random points between $x^*$ and $\hat{x}$, and $\theta^*$ and $\hat{\theta}$ respectively. We introduce $\hat{D_2F} = D_2F(\hat{x},\hat{\theta})$, and an asymptotic expression for $(\theta^*-\hat{\theta})$ is

$$(\theta^*-\hat{\theta})\int_0^1 \hat{D_2F}^\top D_2F(\hat{x},\theta^*) w\, dt = -\int_0^1 \hat{D_2F}^\top (\dot{\hat{x}}-\dot{x}^*) w\, dt$$
$$-\int_0^1 \hat{D_2F}^\top D_1F(\tilde{x}^*,\theta^*)(x^*-\hat{x}) w\, dt.$$

It suffices to consider the convergence in law of the random integral $H_n = \int_0^1 (D_2F^*)^\top \left((\dot{\hat{x}}-\dot{x}^*) + D_1F^*(x^*-\hat{x})\right) w\, dt$ because the random variable

$$G_n = \int_0^1 \hat{D_2F}^\top \left((\dot{\hat{x}}-\dot{x}^*) + D_1F(\tilde{x}^*,\theta^*)(x^*-\hat{x})\right) w\, dt$$

is such that $G_n - H_n \to 0$ in probability (in the $L^2$ sense), moreover we have the convergence in probability of $\int_0^1 \hat{D_2F}^\top D_2F(\hat{x},\theta^*) dt$ to $J^*$, the Hessian of $R_w^2$.

Indeed, we consider the map $\mathcal{D} : (x,\theta) \mapsto (t \mapsto D_2F(x(t),\theta))$ defined on $C([0,1],\mathcal{K}) \times \Theta$ (with the product Hilbert metric) with values in $C([0,1],\mathbb{R}^{n\times p})$ (with the $L^2$ norm $\|A\|_2 = \int_0^1 Tr(A^\top(t)A(t))dt$). Since $D_2F$ is Lipschitz in $(x,\theta)$, the functional map $\mathcal{D}$ is a continuous map, and we can claim by the continuous mapping theorem that the random functions $\hat{D_2F}$ and $t \mapsto D_2F(\hat{x}(t),\theta^*)$ converge in probability (in the $L^2$ sense) to $D_2F^*$. As a consequence, $\|\hat{D_2F}\|_2$ converges in probability to $\|D_2F^*\|_2$ so it is also bounded, and $\|D_2F(\hat{x},\theta^*) - D_2F^*\|_2 \to 0$ in probability. This statement is also true for all entries of these (function) matrices, which enables to claim that all entries of the matrix

$$\int_0^1 \left(\hat{D_2F}\right)^\top \left(D_2F(\hat{x},\theta^*) - D_2F^*\right) w(t) dt$$

tend to zero in probability (by applying the Cauchy-Schwarz inequality componentwise). Moreover, we have convergence in probability of each entry of $\int_0^1 (\hat{D_2F})^\top D_2F^* w\, dt$ to the corresponding entry of $\int_0^1 (D_2F^*)^\top D_2F^* dt$ (consequence of the convergence of $D_2F(\hat{x},\hat{\theta})$ to $D_2F^*$ in the $L^2$ sense), which implies finally that

$$\int_0^1 \left(D_2F(\hat{x},\hat{\theta})\right)^\top D_2F(\hat{x},\theta^*) w\, dt \xrightarrow{P} J^*.$$

By the same arguments and by using the fact that $D_1F$ is also Lipschitz in $(x,\theta)$, we have convergence of the matrix $G_n - H_n$ to 0 in probability. The asymptotic behavior of $(\hat{\theta}_n - \theta^*)$ is then given by the random integral

$$J^{*-1} \int_0^1 (D_2F^*)^\top \left((\dot{\hat{x}}-\dot{x}^*) + D_1F^*(x^*-\hat{x})\right) w\, dt. \tag{3.1}$$

As this functional is linear in $(\dot{\hat{x}} - \dot{x}^*)$ and $(x^* - \hat{x})$, this is equivalent to study the convergence of the functional $\Gamma(\hat{x})$ to $\Gamma(x^*)$, where

$$\Gamma(x) = \int_0^1 (D_2 F^*)^\top \dot{x}(t) w \, dt + \int_0^1 D_2 F^{*\top} D_1 F^* x(t) w \, dt.$$

Since $D_{12}F$ and $\dot{w}$ exist, we can integrate by part the first term of the right-hand side,

$$\int_0^1 (D_2 F^*(t))^\top \dot{x}(t) w(t) dt = \left[ D_2 F^{*\top} x \, w \right]_0^1 - \int_0^1 \frac{d}{dt}(D_2 F^*(t) w(t)) x(t) dt.$$

As a consequence, $\Gamma$ is the sum of the linear functionals $\Gamma_{s,w}$, $\Gamma_{b,w}$ introduced before:

$$\begin{aligned}\Gamma(x) &= \int_0^1 \left( D_2 F^{*\top} D_1 F^* w(t) - \frac{d}{dt}(D_2 F^*(t) w(t)) \right) x(t) dt \\ &\quad + D_2 F^{*\top}(1) x(1) w(1) - D_2 F^{*\top}(0) x(0) w(0).\end{aligned}$$

Finally, by linearity we can write

$$\hat{\theta} - \theta^* = (\Gamma_s(\hat{x}) - \Gamma_s(x^*) + \Gamma_{b,w}(\hat{x}) - \Gamma_{b,w}(x^*)) + o_P(1).$$

□

Proposition 3.1 shows that two-step estimators behave asymptotically as the sum of two plug-in estimators of linear operators: a smooth functional $\Gamma_{s,w}$ and a non-smooth functional $\Gamma_{b,w}$, the latter being a weighted evaluation functional (at the boundaries). It is well-known that the asymptotic behavior of these plug-in estimators are of different nature: Stone [35] showed that the best achievable rate for pointwise estimation over Sobolev classes is $n^{\alpha/2\alpha+1}$ (where $\alpha$ is the number of continuous derivatives of $x$), whereas integral functionals $\int h(t) x(t) dt$ can be estimated at a root-$n$ in a wide variety of situations, e.g. Bickel and Ritov [4]. In particular, we will show in the next section that $\Gamma_{s,w}(x^*)$ can be estimated at a root-$n$ rate with regression splines. A direct consequence of this asymptotic consequence is that two-step estimators does not achieve the optimal parametric rate in generality (with whatever weight function $w$). Hence, two-step estimators used in applications computed with a constant weight function have degraded asymptotic properties. Moreover, the classical bias-variance analysis of pointwise nonparametric estimation enables to give additional insight in the behavior of $\hat{x}$. Indeed, it is well known that the quadratic error $\forall t \in [0,1], E[(\hat{x}(t) - x^*(t))^2] = b(t)^2 + v(t)$, where $b$ is the bias function and $v$ is the variance function. Although $\Gamma_{s,w}(\hat{x}_n)$ is also a biased estimator of $\Gamma_{s,w}(x^*)$, this bias can converge at a faster rate than root-$n$ (under appropriate conditions on $\hat{x}_n$), so that the bias of $\hat{\theta}_n$ is mainly due to the bias terms $b(0)$ and $b(1)$. This remark strongly motivates the use of nonparametric estimators with good boundary properties for reducing the bias of $\hat{\theta}$: a practical consequence of proposition 3.1 is that local polynomials should be preferred to the classical kernel regression estimator because local polynomials have better boundary properties at the boundaries of the interval [13]. Nevertheless, a neat approach to the efficiency of two-step estimators is to restrict to the case of a weight function $w$ with boundary conditions $w(0) = w(1) = 0$, so that $\Gamma_{b,w}$ does vanish and the two step estimator is then asymptotically equivalent to the smooth functional $\Gamma_{s,w}(\hat{x})$, i.e.

$$\hat{\theta} - \theta^* = \Gamma_{s,w}(\hat{x}) - \Gamma_{s,w}(x^*) + o_P(1).$$

Then, it suffices to check that the plug-in estimator $\Gamma_s(\hat{x})$ is a root-$n$ rate estimator of $\Gamma(x^*)$, which depends on the definition of $\hat{x}_n$. We detail then essential properties of regression splines in the next section, and derive the desired plug-in property.

### 3.2. Two-step estimator with least squares splines

A spline is defined as a piecewise polynomial that is smoothly connected at its knots. For a fixed integer $k \geq 2$, we denote $\mathbb{S}(\boldsymbol{\xi}, k)$ the space of spline functions of order $k$ with knots $\boldsymbol{\xi} = (0 = \xi_0 < \xi_1 < \cdots < \xi_{L+1} = 1)$: a function $s$ in $\mathbb{S}(\boldsymbol{\xi}, k)$ is a polynomial of degree $k - 1$, on each interval $[\xi_i, \xi_{i+1}], i = 0, \ldots, L$ and $s$ is in $C^{k-2}$. $\mathbb{S}(\boldsymbol{\xi}, k)$ is a vector space of dimension $L + k$ and the most common basis used in applications are the truncated power basis and the B-spline basis. The latter is usually used in statistical applications as it is a basis of compactly supported functions which are nearly orthogonal. B-splines can be defined recursively from the augmented knot sequence $\boldsymbol{\tau} = (\tau_j, j = 1, \ldots, L + 2k)$ with $\tau_1 = \cdots = \tau_k = 0$, $\tau_{j+k} = \xi_j$, $j = 1, \ldots, L$ and $\tau_{L+k+1} = \cdots = \tau_{L+2k} = 1$. We note $B_{i,k'}$, the $i^{th}$ B-spline basis function of order $k'$ ($1 \leq k' \leq k$) with the corresponding knot sequence $\boldsymbol{\tau}$. The B-spline basis of order $k' = 1, \ldots, k$ are linked then by the recurrence equation:

$$\forall i = 1, \ldots, L + 2k - 1, \forall t \in [0, 1], B_{i,1}(t) = 1_{[\tau_i, \tau_{i+1}]}(t)$$

and $\forall i = 1, \ldots, L + 2k - k', \forall k' = 2, \ldots, k, \forall t \in [0, 1]$,

$$B_{i,k'}(t) = \frac{t - \tau_i}{\tau_{i+k'-1} - \tau_i} B_{i,k'-1}(t) + \frac{\tau_{i+k'} - t}{\tau_{i+k'} - \tau_{i+1}} B_{i+1,k'-1}(t).$$

A (univariate) nonparametric estimator $\hat{x}_n$ is then computed by classical least-squares (the so-called least-squares splines or regression splines). Most of the time, cubic splines ($k = 4$) are used and the essential assumptions are made on the knots sequence $\boldsymbol{\xi}$ in an asymptotic setting: it is supposed that the number of knots $L_n$ tends to infinity with a controlled repartition. The well-posedness of the B-spline basis (and corresponding design matrix) is critical for the good behavior of the regression spline, and the needed material for the asymptotic analysis can be found in [44]. We define below the nonparametric estimator that we use.

We have $n$ observations $y_1, \ldots, y_n$ corresponding to noisy observations of the solution of the ODE (1.1), and we introduce $Q_n$ the empirical distribution of the sampling times times $t_1, \ldots, t_n$. We suppose that this empirical distribution converges to a distribution function $Q$ (which possesses a density $q$ w.r.t Lebesgue measure). We construct our estimator $\hat{x}_n$ of $x^*$ as a function in the space $\mathbb{S}(\boldsymbol{\xi}_n, k)$ of spline functions of degree $k$ and knots sequence $\boldsymbol{\xi}_n = (\xi_0, \xi_1, \ldots, \xi_{L_n+1})$ ($K_n$ is the dimension of $\mathbb{S}(\boldsymbol{\xi}_n, k)$). The knots sequence is chosen such that $\max_{1 \leq i \leq L_n} |h_{i+1} - h_i| = o(L_n^{-1})$, $|\boldsymbol{\xi}|/\min_i h_i \leq M$ where $h_i = (\xi_{i+1} - \xi_i)$ and $|\boldsymbol{\xi}| = \sup_i h_i$ is the mesh size of $\boldsymbol{\xi}$. As a consequence, we have $|\boldsymbol{\xi}| = O(L_n^{-1})$. Like in Zhou *et al.* [44], we suppose that we have convergence of $Q_n$ towards $Q$ at a rate controlled by the mesh size, i.e.

$$\sup_{t \in [0,1]} |Q_n(t) - Q(t)| = o(|\boldsymbol{\xi}|). \tag{3.2}$$

The estimator $\hat{x}_n$ we consider is written componentwise in the basis of B-splines $(B_1, \ldots, B_{K_n})$

$$\forall i = 1, \ldots, d, \ \hat{x}_{n,i} = \sum_{k=1}^{K_n} c_{ik} B_k \tag{3.3}$$

or in matrix form $\hat{x}_n = C_n \mathbf{B}$ with the vector-valued function $\mathbf{B} = (B_1, \ldots, B_{K_n})^\top$ and the $d \times K_n$ coefficient matrix $C_n = (c_{ik}^n)_{1 \leq i, k \leq d, K_n}$ (and column vectors $\mathbf{c}_{i,n} = (c_{i1}, \ldots, c_{i,K_n})^\top \in \mathbb{R}^{K_n}$). We stress the fact that all the components $\hat{x}_{n,i}$ are approximated via the same space, although it may be inappropriate in some practical situations but it enables to keep simple expressions for the estimator. The fact that we look for a function in the vector space spanned by B-splines, puts emphasis on the regression interpretation of the first step of our estimating procedure. The estimation of the parameter $C_n$ can be cast into the classical multivariate regression setting

$$\mathbf{Y}_n = \mathbf{B}_n C_n^\top + \epsilon_n, \tag{3.4}$$

where $\mathbf{Y}_n = (\mathbf{Y}_1 \ldots \mathbf{Y}_d)$ is the $n \times d$ matrix of observations, $\epsilon$ is the $n \times d$ matrix of errors, $C_n^\top$ is the $K_n \times d$ matrix of coefficients and $\mathbf{B}_n = (B_j(t_i))_{1 \leq i \leq n, 1 \leq j \leq K_n}$ is the design matrix. We look for a function close to the data in the $L^2$ sense, i.e. we estimate the coefficient matrix $C_n$ by least-squares

$$\hat{\mathbf{c}}_{i,n} = \arg\min_{\mathbf{c} \in \mathbb{R}^{K_n}} \sum_{j=1}^n (y_{ij} - \mathbf{B}(t_j)^\top \mathbf{c})^2, i = 1, \ldots, d,$$

which gives the least squares estimator $\hat{\mathbf{C}}_n = (\mathbf{B}_n^\top \mathbf{B}_n)^+ \mathbf{B}_n^\top \mathbf{Y}$ where $(\cdot)^+$ denotes the Moore-Penrose inverse. We have

$$\forall i \in \{1, \ldots, d\}, \forall t \in [0,1], \ \hat{x}_{i,n}(t) = \mathbf{B}^\top(t) \hat{\mathbf{c}}_{i,n},$$

where $\hat{\mathbf{c}}_i = (\mathbf{B}_n^\top \mathbf{B}_n)^+ \mathbf{B}_n^\top \mathbf{Y}_i$. Finally, we introduce the projection matrix $P_{B,n} = \mathbf{B}_n (\mathbf{B}_n^\top \mathbf{B}_n)^+ \mathbf{B}_n^\top$. We will use the notation $x \lesssim y$ to denote that there exists a constant $M > 0$ such that $x \leq My$.

General results given by Huang in [21] ensure that $\hat{x}_n \xrightarrow{L^2} x^*$ in probability for sequences of suitably chosen approximating spaces $\mathbb{S}(k, \boldsymbol{\xi}_n)$ with an increasing number of knots. Indeed, corollary 1 in [21] enables us to claim that if the observation times are random with $Q(B) \geq c\lambda(B)$ ($0 < c \leq 1$ and $\lambda(\cdot)$ is the Lebesgue measure on $[0,1]$), the function $x^*$ is in the Besov space $B_{2,\infty}^\alpha$ (with $k \geq \alpha - 1$) and the dimension grows such that $\lim_n \frac{K_n \log K_n}{n} = 0$ then

$$\frac{1}{n}\sum_{i=1}^n (y_{ij} - \hat{x}_i(t_j))^2 + \|\hat{x}_i - x_i^*\|_2 = O_P\left(\frac{K_n}{n} + K_n^{-2\alpha}\right).$$

Moreover, the optimal rate $O_P(n^{-2\alpha/(2\alpha+1)})$ (given by Stone [36]) is reached for $K_n \sim n^{1/(2\alpha+1)}$. For this nonparametric estimator, it is possible to construct a consistent two-step estimator $\hat{\theta}_n$ by minimization of $R_{n,w}^2(\theta)$.

We need then a general result for the asymptotics of linear functionals of spline estimators. This can be seen as a special case of a general result derived by Andrews for series estimators [1]. First, we need to have a precise picture of the evolution of the basis $(B_1, \ldots, B_{K_n})$ as $K_n \to \infty$ and particularly the asymptotic behavior of the empirical covariance $G_{K_n,n} = \frac{1}{n}(\mathbf{B}_n^\top \mathbf{B}_n)$ and of the (theoretical) covariance $G_{K_n} = \int_0^1 \mathbf{B}(t)\mathbf{B}(t)^\top dQ(t)$. From [44], we have that if $K_n = o(n)$,

$$\forall t \in ]0,1], \ \mathbf{B}^\top(t)\left(\mathbf{B}_n^\top \mathbf{B}_n\right)^{-1}\mathbf{B}(t) = \frac{1}{n}\mathbf{B}(t)^\top G_{K_n}^{-1}\mathbf{B}(t) + o\left(\frac{1}{n|\boldsymbol{\xi}|}\right). \quad (3.5)$$

The analysis of section 3.1 gives interest in the asymptotic behavior of $\Gamma_a(\hat{x}_n)$ where $\Gamma_a$ is a linear functional $\Gamma(x) = \int_0^1 a(s)^\top x(s)ds$ where $a$ is a function in $C^1([0,1], \mathbb{R}^d)$. If $x = C^\top \mathbf{B}$, $\Gamma(x) = \sum_{i=1}^d \mathbf{c}_i^\top \boldsymbol{\gamma}_i = Trace(C^\top \boldsymbol{\gamma})$ with $\boldsymbol{\gamma}_i$ the columns of the $K \times d$ matrix $\boldsymbol{\gamma} = \int_0^1 \mathbf{B}(s)A^\top(s)ds$. Hence, the asymptotic behavior is derived directly from the asymptotics of $\hat{\mathbf{C}}_n$ and of matrix $\boldsymbol{\gamma}$. By using the results from Andrews [1], we derive the asymptotic normality of this functional. For simplicity, we consider only the case $d = 1$, the extension to higher dimensions is cumbersome but straightforward. If the variance of the noise is $\sigma^2$, the variance of $\hat{\mathbf{c}}_n$ is

$$V(\hat{\mathbf{c}}_n) = \sigma^2(\mathbf{B}_n^\top \mathbf{B}_n)^+ \quad (3.6)$$

and the variance of the estimator of the functional is

$$V_n = V(\Gamma(\hat{x}_n)) = \sigma^2 \boldsymbol{\gamma}_n^\top (\mathbf{B}_n^\top \mathbf{B}_n)^+ \boldsymbol{\gamma}_n. \quad (3.7)$$

**Proposition 3.2** (**Convergence of linear functionals of regression splines**). *Let $(\boldsymbol{\xi}_n)_{n\geq 1}$ be a sequence of knot sequences of length $L_n + 2$, and $K_n = \dim(\mathbb{S}(\boldsymbol{\xi}_n, k))$, with $k \leq 2$. We suppose that $L_n \to \infty$ is such that $n^{1/2}|\boldsymbol{\xi}_n| \to 0$ and $n|\boldsymbol{\xi}_n| \to \infty$. If $a : [0,1] \to \mathbb{R}$ is in $C^1$, and $x$ is in $C^\alpha$, $2 \leq \alpha \leq k$, $\Gamma(x) = \int_0^1 a(s)x(s)ds$ then:*

**(i)** *$\Gamma(\hat{x}_n) - \Gamma(x) = O_P(n^{-1/2})$ and $\sqrt{n}(\Gamma(\hat{x}_n) - \Gamma(x))$ is asymptotically normal,*

**(ii)** $\forall t \in [0,1]$, $\hat{x}_n(t) - x(t) = O_P(n^{-1/2}|\boldsymbol{\xi}_n|^{-1/2})$,

**(iii)** $V(\hat{x}_n(t))^{-1/2}(\hat{x}_n(t) - x(t))$ *is asymptotically normal,* $t \in [0,1]$.

*Proof.* In order to prove the asymptotic normality of $\Gamma(\hat{x}_n) - \Gamma(x)$, we check the assumptions of theorem 2.1 of [1]. Assumption A is satisfied because the $\epsilon_i$'s are i.i.d. with finite variance. For assumption B, since $a$ is $C^1$, the functional is continuous with respect to the Sobolev norm (or simply the sup norm). From the approximation property of spline spaces [7], it is possible to construct a spline $\tilde{a} = \sum_{i=1}^{K_n} \alpha_i B_i = \boldsymbol{\alpha}_n^\top \mathbf{B} \in \mathbb{S}(\boldsymbol{\xi}_n, k)$ such that $\|a - \tilde{a}\|_\infty = O(|\boldsymbol{\xi}|^1)$. The quality of approximation of the functional $\Gamma_a$ is directly derived: $|\Gamma_a(x) - \Gamma_{\tilde{a}}(x)| = |\int_0^1 (a - \tilde{a})(s)x(s)ds| \lesssim |\boldsymbol{\xi}|\|x\|_\infty$. Hence, it suffices to look at the case $a = \boldsymbol{\alpha}_n^\top \mathbf{B}$ because $\Gamma_a(x) - \Gamma_{\tilde{a}}(x)$ will tend to zero at faster rate than $n^{1/2}$. We introduce the vectors $\gamma_n = (\Gamma_a(B_1) \ldots \Gamma_a(B_{K_n}))^\top$, so we have $\gamma_n^\top \gamma_n = \boldsymbol{\alpha}^\top G_{K_n}^\top G_{K_n} \boldsymbol{\alpha} \geq \lambda_{min}^2 G_{K_n} \times \|\boldsymbol{\alpha}\|_2^2$. From lemma 6.2 of [44] (giving bounds on the eigenvalues of $G_{K_n}$), we get $\gamma_n^\top \gamma_n \gtrsim |\boldsymbol{\xi}|\|\boldsymbol{\alpha}\|_2^2$. Lemma 6.1 from [44] (about the equivalence of $L^2$ norm and Euclidean norm of spline coefficients) ensures that $\gamma_n^\top \gamma_n$ is bounded away from 0 because

$$|\boldsymbol{\xi}|\|\boldsymbol{\alpha}\|_2^2 \gtrsim \int_0^1 a^2(s) dQ_n(s)$$

hence $\liminf_n \gamma_n^\top \gamma_n > 0$ and assumption B is satisfied.

From (3.5), we get the behavior of the diagonal entries of $P_{B,n}$:

$$\forall i \in [1..K_n], \ (P_{B,n})_{ii} = \frac{1}{n}\mathbf{B}(\xi_i)^\top G_{K_n}^{-1} \mathbf{B}(\xi_i) + o\left((n|\boldsymbol{\xi}|)^{-1}\right) \tag{3.8}$$

we see that assumption C(ii) is true because $\mathbf{B}(\xi_i)^\top G_{K_n}^{-1} \mathbf{B}(\xi_i) \leq c_1 \|\mathbf{B}(\xi_i)\|_2^2 |\boldsymbol{\xi}|^{-1}$ and $\|\mathbf{B}(\xi_i)\|_2^2 \leq k$ (because the B-splines are bounded by 1 and only $k$ of them are strictly positive) ensure that $\max_i (P_{B,n})_{ii} = O((n|\boldsymbol{\xi}|)^{-1}) \to 0$. It is clear that $\mathbf{B}_n$ is of full rank for $n$ large enough.

Since $\alpha \leq k$, it exists a sequence $(\tilde{x}_n) \in \mathbb{S}(\boldsymbol{\xi}_n, k)$ such that $\|\tilde{x}_n - x^*\|_\infty = O(|\boldsymbol{\xi}|^\alpha)$ [7], hence

$$n^{1/2}\|\tilde{x}_n - x^*\|_\infty \to 0.$$

If we use again the spline approximation of the function $a$, we derive the following expression for

$$\gamma_n^\top \left(\frac{\mathbf{B}^\top \mathbf{B}}{n}\right)^{-1} \gamma_n = \boldsymbol{\alpha}^\top G_K^\top G_{K,n}^{-1} G_K \boldsymbol{\alpha}.$$

From lemma 6.2. of [44], we have $\boldsymbol{\alpha}^\top G_K^\top G_{K,n}^{-1} G_K \boldsymbol{\alpha} \gtrsim \boldsymbol{\alpha}^\top G_K \boldsymbol{\alpha}$. As for $\gamma_n^\top \gamma_n$, we have

$$\liminf_n \gamma_n^\top \left(\frac{\mathbf{B}^\top \mathbf{B}}{n}\right)^{-1} \gamma_n > 0,$$

which remains true when $a$ is any smooth function in $C^1$.

According to Andrews, we can conclude $V_n^{-1/2}(\Gamma_a(\hat{x}_n) - \Gamma_a(x^*)) \rightsquigarrow N(0,1)$. We obtain an equivalent of the rate of convergence by the same approximation

as above

$$\begin{aligned} V_n &= \gamma_n^\top \left(\mathbf{B}^\top \mathbf{B}\right)^{-1} \gamma_n \\ &= \frac{1}{n}\boldsymbol{\alpha}^\top G_K G_{K,n}^{-1} G_K \boldsymbol{\alpha} \\ &\simeq \frac{1}{n}\boldsymbol{\alpha}^\top G_K \boldsymbol{\alpha} \end{aligned}$$

i.e. $V_n \sim \frac{\|\alpha\|^2 |\boldsymbol{\xi}|}{n}$ and we obtain finally that $V_n \sim n^{-1}$.

The technique used by Andrews for his theorem 2.1 gives also asymptotic normality of $\hat{x}_n(t) = \mathbf{B}(t)^\top \hat{\mathbf{c}}_{i,n}$. We have then

$$\forall t \in [0,1],\ V(\hat{x}_n(t)) = \sigma^2 \mathbf{B}(t)^\top (\mathbf{B}_n^\top \mathbf{B}_n)^+ \mathbf{B}(t)$$

and from (3.5) we get $V(\hat{x}_n(t)) = \frac{\sigma^2}{n}\mathbf{B}(t)^\top G_{K_n}^{-1}\mathbf{B}(t) + o(\frac{1}{n|\boldsymbol{\xi}|})$, so that $V(\hat{x}_n(t)) \sim \frac{C}{n|\boldsymbol{\xi}|}$ from lemma 6.6 in [44]. □

For deriving the asymptotics of the two-step estimator when regression splines are used in the first, we just need to put the results of propositions 3.1 and 3.2 together.

**Theorem 3.1 (Asymptotic normality and rates of the two-step estimator).** *Let $F$ a $C^m$ vector field w.r.t $(\theta, x)$ ($m \geq 1$), such that $D_1 F, D_2 F$ are Lipschitz w.r.t $(\theta, x)$, and $D_{12}F$ exists. We suppose that the Hessian $J^*$ of the asymptotic criterion $R_w^2(\theta)$ evaluated at $\theta^*$ is nonsingular, and that the conditions of proposition 2.1 are satisfied.*

*Let $\hat{x}_n \in \mathbb{S}(\boldsymbol{\xi}_n, k)$ a regression spline with $k \geq 2$, such that $n^{1/2}|\boldsymbol{\xi}_n| \to 0$ and $n|\boldsymbol{\xi}_n| \to \infty$, then the two-step estimator $\hat{\theta}_n = \arg\min_\theta R_{n,w}^2(\theta)$ is asymptotically normal and*

**(i)** *if $w(0) = w(1) = 0$, then $(\hat{\theta}_n - \theta^*) = O_P(n^{-1/2})$*
**(ii)** *if $w(0) \neq 0$ or $w(1) \neq 0$, then $(\hat{\theta}_n - \theta^*) = O_P(n^{-1/2}|\boldsymbol{\xi}_n|^{-1/2})$.*

*In the case (ii), the optimal rate of convergence for the Mean Square Error is obtained for $K_n = O(n^{1/(2m+3)})$ and we have then $(\hat{\theta}_n - \theta^*) = O_P(n^{-(m+1)/(2m+3)})$.*

*Proof.* From proposition 3.1, we have

$$\hat{\theta} - \theta^* = J^{*-1}\left(\Gamma_{s,w}(\hat{x}) - \Gamma_{s,w}(x^*) + \Gamma_{b,w}(\hat{x}) - \Gamma_{b,w}(x^*)\right) + o_P(1).$$

so that we just have to apply proposition 3.2 to $\Gamma_{s,w}(\hat{x})$ and $\Gamma_{b,w}(\hat{x})$. We can claim the asymptotic normality of $\sqrt{n}(\Gamma_s(\hat{x}_n) - \Gamma_s(x^*))$ and of $\sqrt{n|\boldsymbol{\xi}_n|}(\Gamma_{b,w}(\hat{x}_n) - \Gamma_{b,w}(x^*))$ (normality is extended from scalar functional to multidimensional functional by the Cramér-Wold device). We have then two cases for the rate of convergence of $\hat{\theta}_n$, depending on the values $w(0), w(1)$. When $w(0) = w(1) = 0$, there is only the parametric part, but when $\Gamma_{b,w}$ does not vanish the nonparametric part with rate $\sqrt{n|\boldsymbol{\xi}_n|}$ remains. We can determine the optimal rate of convergence in the mean square sense by using the Bias - Variance decomposition for the evaluation functional $\|\hat{\theta}_n - \theta^*\|^2 = O_P\left((E(\hat{x}_n(t)) - x(t))^2\right) +$ $O_P(Var(\hat{x}_n(t)))$. Theorem 2.1 of [44] gives $E\left((\hat{x}_n(t) - x^*(t))^2\right) = O(|\boldsymbol{\xi}_n|^{m+1})$ (because $x^*$ is $C^{m+1}$) and $Var(\hat{x}_n(t)) = O_P(n^{-1}|\boldsymbol{\xi}_n|^{-1})$ so the optimal rate is reached for $|\boldsymbol{\xi}_n| = O(n^{-1/(2m+3)})$ and is $O(n^{-(2m+2)/(2m+3)})$. □

**Remark 3.1** (**Random observational times**). The asymptotic result given for the deterministic observational times $0 \leq t_1 < \cdots < ct_n \leq 1$ remains true when they are replaced by realizations of some random variables $T_1, \ldots, T_n$ as long as the assumptions of the two previous propositions are true with probability one. Andrews gives some conditions (theorem 2) in order to obtain this. It turns out that in the case of $T_1, \ldots, T_n$ i.i.d. random variables drawn from the distribution $Q$, it suffices to have $K_n^4 \lesssim n^r$ with $0 < r < 1$. In particular, as soon as $m \geq 1$, the conclusion of proposition 3.1 holds with probability one for the optimal rate $K_n = n^{1/(2m+3)}$.

We have restricted theorem 3.1 to regression splines in order to a have precise and self-content statement of conditions under which the asymptotics of the two-step estimator is known. In particular cubic splines gives root-$n$ consistent estimators for smooth differential equation ($m \geq 2$) and appropriate weight function.

The main point of our study is that the asymptotic normality and parametric rate are derived from the behavior of the smooth functional $\Gamma_{s,w}$, which can be derived for series estimators (polynomials or Fourier series) from [1]. Moreover, the same theorem can be derived for kernel regression (Nadaraya-Watson) by using the results on plug-in estimators of Goldstein and Messer [18]. More generally, the same theorem may be obtained for other nonparametric estimators such as local polynomial regression, orthonormal series, or wavelets.

## 4. Experiments

The Lotka-Volterra equation is a standard model for the evolution of prey-predator populations. It is a planar ODE

$$\begin{cases} \dot{x} &= ax + bxy \\ \dot{y} &= cy + dxy \end{cases} \quad (4.1)$$

whose behavior is well-known [19]. Despite its simplicity, it can exhibit convergence to limit cycles which is one of the main features of nonlinear dynamical systems, which has usually a meaningful interpretation. Due to this simplicity and the interpretability of the solution, it is often used in biology (population dynamics or phenomena with competing species), but the statistical estimation of parameters $a$, $b$, $c$, $d$ in 4.1 has not been extensively addressed. Nevertheless, Varah (1982) illustrates the two-step method with cubic splines and knots chosen by an expert on the same model as (4.1). Froda *et al.* (2005) [16] have considered another original method exploiting the fact that the orbit may be a closed curve for some values of the parameters. In this section, we will consider a slight generalization and reparametrization of the previous model which consists in adding the squared terms $x^2$ and $y^2$:

$$\begin{cases} \dot{x} &= x(a_1 x + a_2 y + a_3) \\ \dot{y} &= y(b_1 x + b_2 y + b_3) \end{cases} \tag{4.2}$$

We use the two-step estimator obtained by minimizing $R_{n,w}^2(\theta)$ in order to illustrate the consistency and the asymptotic normality of the estimator proved in the previous section. In particular, we will consider two estimators: one obtained with a uniform weight function $w = 1$, and a second one with a weight vanishing at the boundaries. According to theorem 3.1, it gives two different rates of convergence: we consider then the influence of the number of observations $n = 20, 30, 50, 100, 200, 500, 1000$. We consider also a small number of observations ($n = 20, 30$) as the nonparametric procedure can give poor results in the small-sample case, and the simulation can give an insight into the expected properties in this practical situation. As the reconstruction of the solution in the first step is critical, we consider the estimation of an ODE with two different parameters $\theta_1$ with $a_1 = 0$, $a_2 = -1.5$, $a_3 = 1$, $b_1 = 2$, $b_2 = 0$ and $b_3 = -1.5$, and $\theta_2$ with $a_1 = 0$, $a_2 = -1.5$, $a_3 = 1$, $b_1 = 2$, $b_2 = 1$ and $b_3 = -1.5$. In both cases, the parameters of the quadratic terms $a_1$ and $b_2$ are supposed to be known, so that 4 parameters have to be estimated from noisy discrete observations of the solution of the ODE. We suppose also that the initial conditions are not known, so that they are nuisance parameters that are not estimated by our procedure. These two parameters $\theta_1, \theta_2$ gives rise to two different qualitative behavior of the solution as it can be seen in figures 1, 2 and it gives a view of the influence of their shapes on the inference. The shape of the solution has two consequences: the identifiability of the model through the asymptotic criterion $R_w^2$, and the difficulty of the first step (if the curve is rather wiggly or flat). The data are simulated according to 2.1 where $\epsilon$ is a Gaussian white noise with standard-deviation $\sigma = 0.2$, and the observation times are deterministic and equal to $t_j = j\frac{20}{n}$, $j = 0, \ldots, n-1$.

We define now exactly the two-step estimator that have been used in the experiments. In the first step, we have used a least square cubic spline with an increasing number of knots $K_n$. The selection of $K_n$ and of the position of the knots is a classical problem in the spline regression literature. The study of the asymptotics gives only the rate of $K_n$ (or $|\boldsymbol{\xi}_n|$) in order to optimize the rate of convergence, but one has to find practical devices in order to compute a correct estimator. For instance, this choice was left to the practitioner in Varah's paper, but due to the intensive simulation framework for the estimation of the bias and variance of two-step estimators, we have used an adaptive selection of the knots in order to obtain an automatic first step with reliable and stable approximation properties. This procedure permits then to focus on the quality of the two-step estimator, and not in particular on the first step. We have considered 3 different basis of splines with uniformly spaced knots $\xi_j = j\frac{20}{K_n + 1}$, $j = 0, \ldots, K_n$:

- if $n = 20, 30$, $K_n = 15$,
- if $n = 50$, $K_n = 20$,
- if $n = 100, 200, 500, 1000$, $K_n = 30$.

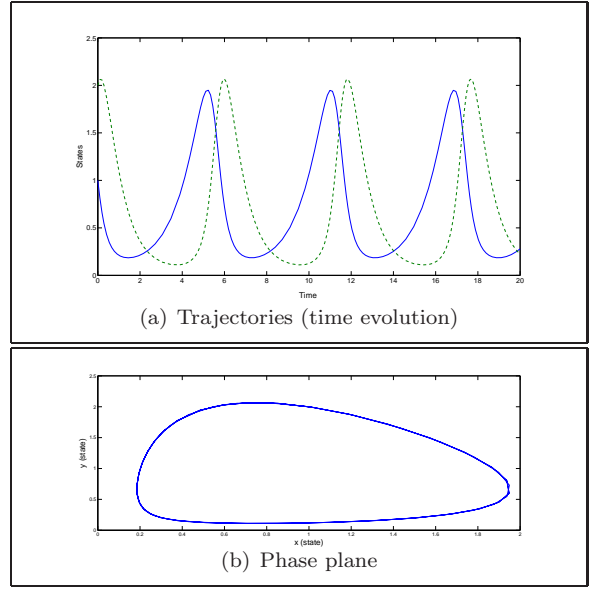

FIGURE 1. *Solution of Lotka-Volterra system in the phase plane for $\theta_1$ and $x(0) = 1$ et $y(0) = 2$.*

and we have selected the knots by applying a variable selection approach to the truncated power basis (instead of the B-spline basis) defined as the set of functions $1, t, \ldots, t^{k-1}, (t - \xi_1)_+^{k-1}, \ldots, (t - \xi_L)_+^{k-1}$. The knots are selected by minimizing the Generalized Cross Validation criterion (GCV):

$$GCV(\mathcal{K}_m) = \sum_{i=1}^{n} \frac{y_i - \phi(t_i)}{1 - d(\mathcal{K}_m)/n}$$

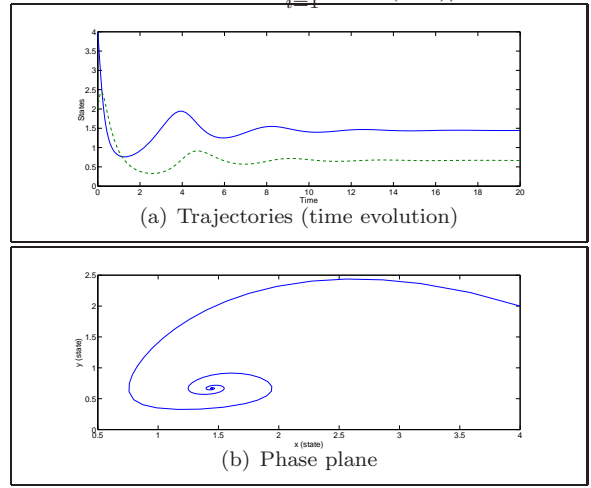

FIGURE 2. *Solution of Lotka-Volterra system in the phase plane for $\theta_2$ and $x(0) = 4$ et $y(0) = 2$.*

TABLE 1
Case 1: Mean and standard deviation of the two-step estimator with a weight function
vanishing at the boundaries ($\hat{\theta}_{1,w}$) and with uniform weigth ($\hat{\theta}_1$). True parameter is
$\theta_1 = (-1.5, 1, 2, -1.5)$

| $n$ | mean($\hat{\theta}_{1,w}$) | mean($\hat{\theta}_1$) | std($\hat{\theta}_{1,w}$) | std($\hat{\theta}_1$) |
|---|---|---|---|---|
| 20 | (−0.97, 0.68, 0.99, −0.75) | (−0.87, 0.57, 0.95, −0.77) | (0.13, 0.12, 0.17, 0.17) | (013, 0.12, 0.17, 0.17) |
| 30 | (−0.95, 0.71, 1.12, −0.95) | (−0.97, 0.72, 1.10, −0.97) | (0.10, 0.10, 0.13, 0.12) | (0.11, 0.10, 0.13, 0.12) |
| 50 | (−1.12, 0.81, 1.41, −1.14) | (−1.11, 0.81, 1.38, −1.16) | (0.11, 0.10, 0.14, 0.13) | (0.11, 0.10, 0.14, 0.13) |
| 100 | (−1.28, 0.88, 1.73, −1.34) | (−1.26, 0.87, 1.71, −1.36) | (0.14, 0.11, 0.18, 0.16) | (0.13, 0.10, 0.19, 0.15) |
| 200 | (−1.37, 0.92, 1.83, −1.4) | (−1.35, 0.91, 1.81, −1.41) | (0.09, 0.07, 0.12, 0.1) | (0.10, 0.07, 0.13, 0.1) |
| 500 | (−1.42, 0.95, 1.90, −1.44) | (−1.42, 0.95, 1.90, −1.45) | (0.06, 0.05, 0.08, 0.06) | (0.07, 0.05, 0.09, 0.06) |
| 1000 | (−1.45, 0.97, 1.94, −1.46) | (−1.45, 0.96, 1.93, −1.46) | (0.04, 0.03, 0.06, 0.05) | (0.045, 0.03, 0.07, 0.05) |

TABLE 2
Case 1: Mean Squared Errors of the weighted and unweighted parametric estimators versus
the Mean Squared Error of the nonparametric estimators of the 2 curves

| $n$ | RMSE($\hat{\theta}_{1,w}$) | RMSE($\hat{\theta}_1$) | RMSE($\hat{\phi}$) |
|---|---|---|---|
| 20 | 1.41 | 1.5 | (1.32, 1.56) |
| 30 | 1.21 | 1.21 | (1.05, 1.21) |
| 50 | 0.83 | 0.84 | (0.83, 0.9) |
| 100 | 0.47 | 0.49 | (0.54, 0.55) |
| 200 | 0.3 | 0.32 | (0.38, 0.4) |
| 500 | 0.18 | 0.19 | (0.25, 0.26) |
| 1000 | 0.12 | 0.13 | (0.18, 0.20) |

where $\mathcal{K}_m$ is a subset (of size $m$) of $\{\xi_j, j = 0, \ldots, K_n\}$, and $d(\mathcal{K}_m)$ is the effective number of parameters. As in [15, 27], we use $d(\mathcal{K}_m) = 3m + 1$, and we use an heuristic for considering efficiently the various subsets, which is based on the elimination of the less informative knots (for GCV), and the addition of the more informative knots in order to minimize the GCV criterion (*ElimAdd* procedure in [27]). This procedure is not optimal but it gives a simple and reliable adaptive nonparametric estimators. Other knots selection procedure would have been based on free-knot splines [8, 37]. Experiments have been also performed in the non-adaptive case [5].

The quality of the nonparametric estimation procedure is measured by the

TABLE 3
Case 1: Minima of the criteria $R_n$

| $n$ | $R_{n,w}^2(\hat{\theta}_{1,w})$ | $R_n^2(\hat{\theta}_1)$ |
|---|---|---|
| 20 | (6.86, 9.35) | (11.4, 15.4) |
| 30 | (3.44, 4.11) | (4.25, 5.25) |
| 50 | (2.41, 2.97) | (2.94, 3.87) |
| 100 | (1.86, 2.46) | (2.31, 3.67) |
| 200 | (1.25, 1.62) | (1.62, 2.68) |
| 500 | (0.70, 0.96) | (0.96, 1.56) |
| 1000 | (0.47, 0.64) | (0.64, 1.0) |

TABLE 4
Case 2:Bias and standard deviation of two-step estimators with a weight function vanishing at the boundaries ($\hat{\theta}_{2,w}$) and with uniform weight ($\hat{\theta}_2$). True parameter is
$\theta_2 = (-1.5, 1, 1.5, -1.5)$

| $n$ | mean$\hat{\theta}_{2,w}$ | mean$\hat{\theta}_2$ | std($\hat{\theta}_{2,w}$) | std($\hat{\theta}_2$) |
|---|---|---|---|---|
| 20 | (−0.88, 0.56, 0.48, 0.03) | (−1.18, 0.74, 0.22, 0.38) | (0.27, 0.18, 0.36, 0.52) | (0.28, 0.20, 0.35, 0.49) |
| 30 | (−1.18, 0.77, 0.72, −0.34) | (−1.42, 0.92, 0.46, 0.04) | (0.25, 0.18, 0.24, 0.35) | (0.22, 0.16, 0.19, 0.26) |
| 50 | (−1.17, 0.78, 1.15, −0.94) | (−1.58, 1.03, 0.99, −0.69) | (0.18, 0.13, 0.33, 0.49) | (0.29, 0.20, 0.44, 0.64) |
| 100 | (−1.39, 0.92, 1.12, −0.93) | (−1.45, 0.95, 1.01, −0.78) | (0.18, 0.12, 0.33, 0.48) | (0.18, 0.19, 0.53, 0.77) |
| 200 | (−1.43, 0.95, 1.29, −1.19) | (−1.54, 1.0, 1.17, −1.00) | (0.13, 0.09, 0.25, 0.37) | (0.24, 0.16, 0.43, 0.63) |
| 500 | (−1.45, 0.97, 1.43, −1.40) | (−1.58, 1.05, 1.34, −1.26) | (0.08, 0.05, 0.16, 0.24) | (0.19, 0.12, 0.31, 0.45) |
| 1000 | (−1.47, 0.98, 1.47, −1.45) | (−1.58, 1.05, 1.41, −1.37) | (0.05, 0.04, 0.10, 0.14) | (0.14, 0.09, 0.20, 0.29) |

Root Mean Squared Error (RMSE)

$$RMSE = \left( \int_0^{20} \left( \hat{x}_j - x_j^* \right)^2 (t) dt \right)^{1/2}$$

and is given in tables 2, 5, so that the MSE of the two-step estimator can be related directly to the quality of the first step.

In a general manner, the second step can be addressed by using whatever (nonlinear) optimization procedures. Nevertheless, in the case of the Lotka-Volterra model, it turns out that we have a closed-form expression for $\hat{\theta}$. Indeed, the criterion $R^2_{n,w}$ for the first dimension can be written

$$\begin{aligned} R^2_{n,w}(\theta) &= \int_0^{20} (\dot{\hat{x}} - \hat{x}(a_1\hat{x} + a_2\hat{y} + a_3))^2 w(t) dt \\ &\approx \sum_j \Delta_j w(t_j) \hat{x}(t_j) \left( \frac{\dot{\hat{x}}}{\hat{x}}(t_j) - a_1\hat{x}(t_j) + a_2\hat{y}(t_j) + a_3 \right)^2 \\ &= \sum_j \tilde{w}_j \left( \frac{\dot{\hat{x}}}{\hat{x}}(t_j) - (a_1\hat{x}(t_j) + a_2\hat{y}(t_j) + a_3) \right)^2 \end{aligned}$$

where the criterion is approximated by a Riemann sum ($\Delta_j$ is the size of the subdivision, supposed to be uniform) and $\tilde{w}_j = \Delta_j w(s_j) \hat{x}(s_j)$ is a new weight function. Hence, the estimator is the solution of a weighted least-square program whose solution exists, is unique and has a known expression. This means that the two-step estimator remarkably furnishes a fast and reliable procedure, where there is no problem of local minima, although we deal with a nonlinear regression model.

The quality of the two-steps estimator has been evaluated by computing its mean and standard deviation through a Monte-Carlo study with 1000 independent drawings. The results are shown in tables 1, 4 for the models with $\theta_1$ and $\theta_2$, and illustrates the consistency of the two-step procedures as $n$ grows to infinity. We compute with the same nonparametric estimator $\hat{x}$, and either with a uniform weight function $w = 1$, or boundaries vanishing weight. In our experiment, we use a piecewise linear function with $w(0) = 0$, $w(1) = 1$, $w(19) = 1$ and $w(20) = 0$ (and $w(t) = 1$, $t \in [1, 19]$). We have also a picture of the behavior of the estimator for small sample size and it appears that the two-step estimator (weighted or unweighted) is biased, for cases 1 and 2, and the bias diminishes significantly when $n \leq 100$. An important feature is that the weighted estimator is better behaved than the unweighted one, in both experiments and for small $n$ and big $n$. Indeed, the standard deviation of $\hat{\theta}_{1,w}$ is equal (in case 1) or smaller (in case 2) to the one of $\hat{\theta}_1$, but the main difference between the two estimators comes from the bias term which induces a bigger RMSE for the unweighted estimator, see tables 2, 5. This difference comes from the presence of bias term in the evaluation functionals in the unweighted estimator (as it is emphasized in theorem 3.1), which can be particularly important at the boundaries. In case 2, this difference is very important, due to the flat part for $t \geq 10$.

TABLE 5
*Case 2: Mean Squared Errors of the weighted and unweighted parametrics estimators versus the mean squared of the nonparametric estimators of the 2 curves*

| $n$ | MSE($\hat{\theta}_{2,w}$) | MSE($\hat{\theta}_2$) | MSE($\hat{\phi}$) |
|---|---|---|---|
| 20 | 2.0 | 2.34 | (0.93, 0.97) |
| 30 | 1.48 | 1.88 | (0.68, 0.73) |
| 50 | 0.91 | 1.14 | (0.57, 0.59) |
| 100 | 0.83 | 1.18 | (0.44, 0.40) |
| 200 | 0.52 | 0.87 | (0.53, 0.34) |
| 500 | 0.28 | 0.56 | (0.20, 0.20) |
| 1000 | 0.17 | 0.37 | (0.15, 0.14) |

TABLE 6
*Case 2: Minima of the criterions*

| $n$ | $R_{n,w}^2(\hat{\theta}_{2,w})$ | $R_n^2(\hat{\theta}_2)$ |
|---|---|---|
| 20 | (3.66, 1.99) | (5.8, 2.80) |
| 30 | (2.35, 1.16) | (3.4, 1.57) |
| 50 | (2.01, 1.55) | (3.46, 2.29) |
| 100 | (0.89, 0.40) | (1.60, 0.69) |
| 200 | (0.52, 0.34) | (1.20, 0.61) |
| 500 | (0.30, 0.19) | (0.95, 0.36) |
| 1000 | (0.18, 0.11) | (0.66, 0.21) |

From the expression of the $R_{n,w}^2$, it is clear that the quality of the $\hat{\theta}$ is directly related to the distance $\|\hat{x} - x\|_2^2$, as one can see in tables 2, 5 but our experiment shows that it is not sufficient. First, one need also to estimate correctly the derivative of the solution, and second a low minimum of the criterion ($R_{n,w}^2(\hat{\theta}_{n,w})$) does not indicate a good estimator. This difficulty arises when we compare the distance of the nonparametric estimator and the value of the criterion (see tables 3, 6). For instance, we have a better nonparametric estimator and criterion in case 2 than in case 1 (e.g. $n = 100$), but we have a lower RMSE in case 1. Hence, case 2 appears as a more difficult model to estimate, and the shape of the solution has a direct influence on the asymptotic criterion $R_w^2$ and on the ability to approximate it. We have also controlled the normality of the two-step estimator with a univariate Kolmogorov-Smirnov test for the four parameters. The test used was adapted to the case of unknown mean and variance (via Lilliefors procedure), and the normality assumption cannot be rejected (at a level of 5%) as soon as $n = 20$ for all parameters.

## 5. Conclusion

We have proposed a new family of parametric estimators of ODE's relying on nonparametric estimators, which are simpler to compute than straightforward parametric estimators such as MLE or LSE. The construction of this parametric estimator puts emphasis on the regression interpretation of the ODE's estimation problem, and on the link between a parameter of the ODE and an associated function. By using an intermediate functional proxy, we expect to gain information and precision on likely value of the parameters. We do not have studied the effect of using shape or inequality constraints of the estimator $\hat{x}_n$ but it might be valuable information for the inference of complex models, either by shortening the computation time (it gives more suitable initial conditions) or by accelerating the rate of convergence of the estimator thanks to restriction to smaller sets of admissible parameter values.

We have particularly studied the case $R^2_{n,w}(\theta)$, but other M-estimators such as the one obtained from $R^1_{n,w}(\theta)$ may possess interesting theoretical and practical properties such as robustness. This could be particularly useful in the case of noisy data which can give oscillating estimates of the derivatives of the function.

We have given a general study of the two-step estimator (consistency, asymptotic expansion), and we have shown that the weight function used in $R^2_{n,w}$ controls dramatically the rate of convergence of the estimator, and that this method can furnish root-$n$ consistent estimators. We have then provided a detailed account of the asymptotic behavior of spline-based estimators. This choice is mainly due to the practical interest and the wide use of splines, but the conclusions may remain the same for usual nonparametric estimators (series estimators, kernel estimators). Indeed, we have shown that the asymptotic behavior of the two-step estimator comes from the behavior of a linear functional of the regression function.

In the experiments, we have illustrated the influence of the weight function, as the influence of the solution in the quality of the estimators. In particular, we have shown that the approximation quality of the solution is not sufficient in order to have a good estimator, and that it depends on the shape of the solution.

**Acknowledgement**

The author is grateful to Florence d'Alché-Buc (Université d'Evry) for defining this postdoctoral research project and for stimulating discussions all along this work, and to Chris Klassen (Universiteit van Amsterdam) for his advice about the proofs. Their comments greatly improved the quality of the paper.